\documentstyle{article}
\newtheorem{Definition}{Definition }[section]
\newtheorem{Theorem}[Definition]{Theorem }
\newtheorem{Lemma}[Definition]{Lemma }
\newtheorem{Proposition}[Definition]{Proposition }

\def\proof{\par \noindent{\bf Proof. }\nopagebreak}

\newcommand{\Remark}[1]{\vspace{.2cm}\par\noindent{\bf Remark }{#1}
\vspace{.2cm}}

\def\example{\vspace{.2cm}\par \noindent{\bf Example }}
\newcommand{\claim}[1]{\vspace{.2cm}\par \noindent{\sc Claim }{#1}
\vspace{.2cm}}
\def\eproof{\nopagebreak\par\nopagebreak\hfill {\Large $\diamond$}\par}

\def\O{{\cal O}}

\def\Proj{{\bf P}}
\def\O{{\cal O}}
\def\Q{{\bf Q}}

\def\N{{\bf N}}

\def\C{{\bf C}}

\def\f{\varphi}
\def\ra{\rightarrow}

  \def\flip{\stackrel{\scriptstyle - - >}{  }}

\title{Existence of good divisors on Mukai varieties}
\author{Massimiliano Mella}
\date{}
\begin{document}
\maketitle
\section*{Introduction}
A normal projective variety $X$ is called {\sf Fano} if a multiple
of the \hbox{anticanonical} Weil divisor, $-K_X$, is an ample Cartier
divisor.
The importance of Fano \hbox{varieties} is twofold, from one side they give,
has predicted by Fano \cite{Fa}, \hbox{examples} of non rational varieties
having plurigenera and irregularity all zero (cfr  \cite{Is});
on the other hand they should be the building block of
uniruled variety. Indeed recently, Minimal Model Theory predicted that
every
uniruled variety is birational to a fiber space whose general fiber is
a Fano variety with terminal singularities (cfr \cite{KMM}).

The index of a Fano variety $X$ is the number
$$i(X):=sup\{t\in \Q:
-K_X\equiv tH,\mbox{\rm     for some ample Cartier divisor $H \}$}.$$ It is
known that
$0<i(X)\leq dimX+1$ and if $i(X)\geq dim X$ then $X$ is either an
hyperquadric or a projective space by the Kobayashi--Ochiai criterion.
Smooth Fano n-folds
of index $i(X)=n-1$, {\sf del Pezzo n-folds},
have been classified by Fujita \cite{Fu} and terminal Fano n-folds of
index $i(X)>n-2$ have been independently classified by
Campana--Flenner \cite{CF} and Sano \cite{Sa}.

If $X$ has log terminal singularities,
then $Pic(X)$ is torsion free and therefore, the $H$ satisfying
$-K_X\equiv i(X)H$ is uniquely determined and is called the {\sf
fundamental divisor} of $X$.
Mukai announced, \cite{Mu}, the classification of smooth Fano n-folds
$X$
of index $i(X)=n-2$, under the assumption that the linear system $|H|$
contains a smooth divisor.

It is usually said that a
Fano variety $X$ has {\sf good divisors} if the generic element
of the fundamental divisor of $X$ has at worst the same singularities
as $X$.

Our main Theorem is the following

\vspace{.5cm}\noindent {\bf Theorem 1} {\it Let $X$ be a
Mukai
variety with at worst log terminal singularities.
Then $X$ has good divisors except in the following cases:
\begin{itemize}
\item[-]  $X$ is a singular terminal Gorenstein 3-fold which is
 a complete intersection of a quadric
and a sestic in $\Proj(1,1,1,1,2,3)$
\item[-] $X$ is a terminal not Gorenstein 3-fold
 and the canonical cover of $X$ is a complete
intersection of a quadric and a quartic in $\Proj(1,1,1,1,1,2)$.
\end{itemize}
In both exceptional cases the generic element of the
fundamental divisor has canonical singularities.} \vspace{.5cm}

In particular this proves Mukai hypothesis and
therefore the result of Mukai \cite{Mu} provide a complete
classification of smooth Fano n-folds of index $i(X)=n-2$, {\sf Mukai
manifolds}, see
also \cite{CLM} for a different approach.

The ancestors of the theorem, and indeed the lighthouses that
directed
its proof, are Shokurov's proof
for smooth Fano 3-folds, \cite{Sh} and Reid's extension to
canonical Gorenstein 3-folds using the Kawamata's base point free
technique \cite{Re}. This technique was then applied by Wilson in the
case of smooth Fano
\hbox{4-folds} of index 2, \cite{Wi}, afterwards Alexeev, \cite{Al} did it
for
log terminal Fano \hbox{n-folds} of index $i(X)>n-2$ and  recently Prokhorov
used it to prove Mukai Conjecture in dimension 4 and 5, \cite{Pr1}
\cite{Pr2} \cite{Pr3}. Theorem 1 was proved for smooth
Mukai variety in a first version of this paper,\cite{Me}, 
using Helmke's inductive procedure,\cite{He}. Our main tools will be Kawamata's
base point free technique,
Kawamata's
notion of centers of log canonical singularities, \cite{Ka1} and his
subadjunction formula for codimension $\leq 2$ minimal centers \cite{Ka2}.
These tools allows to replace difficult
non vanishing arguments by a simple Riemann--Roch calculation.

While working on this subject I had several discussions with M.
Andreatta, I would like to express him my deep gratitude. I would
also like to thank A. Corti and Y. Prokhorov
for valuable comments and Y. Kawamata for signaling a gap in the first version
of this paper.
This research was partially supported by  the Istituto Nazionale di
Alta Matematica Francesco Severi  (senior grant 96/97) and
 the Consiglio Nazionale delle Ricerche.

\section{Preliminary results}

We use the standard notation from algebraic geometry.
In particular it is compatible with that of \cite{KMM}
to which we refer constantly, everything is defined over \C.

In the following $\equiv$ (respectively $\sim$) will
indicate
numerical  (respectively linear) equivalence of divisors.
Let $\mu:Y\ra X$ a birational morphism of normal varieties. If
$D$ is a \Q-divisor (\Q-Cartier) then is well defined the strict
transform $\mu_*^{-1}D$ (the pull back $\mu^*D$). For a pair $(X,D)$
of a variety
$X$ and a \Q-divisor $D$, a log resolution is a proper birational
morphism $\mu:Y\ra X$ from a smooth $Y$ such that the union of the
support of $\mu_*^{-1}D$ and of the exceptional locus is a normal
crossing divisor.
\begin{Definition} Let $X$ be a normal variety and $D=\sum_id_iD_i$
an effective
\Q-divisor such that $K_X+D$ is \Q-Cartier. If $\mu:Y\ra X$ is a
log resolution of the pair $(X,D)$, then we can write
$$K_Y+\mu_*^{-1}D=\mu^*(K_X+D)+F$$
with $F=\sum_je_jE_j$ for the exceptional divisors $E_j$. We call
$e_j\in \Q$ the discrepancy coefficient for $E_j$, and regard $-d_i$
as the discrepancy coefficient for $D_i$.

The pair $(X,D)$ is said to have {\sf log canonical} (LC)
(respectively {\sf Kawamata log terminal} (KLT), {\sf log terminal}
(LT)) singularities if
$d_i\leq 1$ (resp. $d_i< 1$, $d_i=0$) and $e_j\geq -1$ (resp.
$e_j>-1$, $e_j>-1$) for any $i,j$ of a log resolution $\mu:Y\ra
X$. The {\sf
log canonical threshold} of a pair $(X,D)$ is $lct(X,D):=
sup\{t\in\Q$: $(X,tD)$ is LC$\}$.
\label{lc}
\end{Definition}
\begin{Definition} A {\sf log-Fano variety} is a pair $(X,\Delta)$
with KLT singularities and such that for some positive integer $m$,
$-m(K_X+\Delta)$ is an ample Cartier divisor. The index of a log-Fano
variety $i(X,\Delta):=sup \{t\in \Q: -(K_X+\Delta)\equiv tH$ for some
ample Cartier divisor $H \}$ and the $H$ satisfying
$-(K_X+\Delta)\equiv i(X,\Delta)H$ is called fundamental divisor.
In case $\Delta=0$ we have log terminal Fano variety.
\end{Definition}

\begin{Proposition}[\cite{Al}] Let $(X,\Delta)$ be a log-Fano n-fold of index
r, $H$ the fundamental divisor in $X$. If $r\geq n-2$ then $dim |H|\geq n-1$.
\label{al}
\end{Proposition}

Let us recall the notion and properties of minimal
centers of log canonical singularities as introduced in \cite{Ka1}
\begin{Definition}[\cite{Ka1}] Let $X$ be a normal variety and $D=\sum
d_iD_i$ an effective \Q-divisor such that $K_X+D$ is \Q-Cartier.
A subvariety $W$ of $X$ is said to be a {\sf center of log canonical
singularities} for the pair $(X,D)$, if there is a birational morphism
from a normal variety $\mu:Y\ra X$ and a prime divisor $E$ on $Y$ with
the discrepancy coefficient $e\leq -1$ such that $\mu(E)=W$.
The set of all centers of log canonical singularities is denoted
by $CLC(X,D)$.
\end{Definition}
\begin{Theorem}[\cite{Ka1},\cite{Ka2}] Let $X$ be a normal variety and $D$ an
effective \Q-Cartier divisor such that $K_X+D$ is \Q-Cartier. Assume
that $X$ is LT and $(X,D)$ is LC.
\begin{itemize}
\item[i)]
If $W_1,W_2\in CLC(X,D)$ and $W$ is
an irreducible component of $W_1\cap W_2$, then $W\in CLC(X,D)$. In
particular, if $(X,D)$ is not KLT then there
exist minimal elements in $CLC(X,D)$.
\item[ii)]
If $W\in CLC(X,D)$ is a minimal center then $W$ is normal
\item[iv)] (subadjunction formula) If $W$ is a minimal center
for $CLC(X,D)$
and $cod W\leq 2$ then there exists an effective
\Q-divisors $\Delta$ on $W$ such that $(K_X+D)_{|W}\equiv
K_W+\Delta$ and
$(W,\Delta)$ is KLT.
\end{itemize}
\label{clc}
\end{Theorem}

\section{Proof of Theorem 1}
Let me sketch the idea of the proof, to make it more
transparent.
By Bertini Theorem if $X$ has
not good divisors then the generic element of $|H|$ has a center of
"bad" singularities contained in $Bsl|H|$. We will derive a
contradiction producing a section of $|H|$ not vanishing identically
on $Z$. The following lemma will be frequently used to this purpose.

\begin{Lemma} Let $X$ be a log terminal Fano n-fold, with $n\geq 3$
and $H$ an ample Cartier divisor with $-K_X\equiv
(n-2)H$ and $G$ a \Q-Cartier divisor. Assume that $(X,G)$ is LC,
$Z\in CLC(X,G)$ is a minimal center  and $G\equiv
\gamma H$, with $\gamma< cod Z-1$.
Then there is a section
of $H$ not vanishing identically on $Z$.
\label{KV}
\end{Lemma}
\proof
First let us perturb $G$ to construct a
\Q-divisor $G_1\equiv \gamma_1 H$ such that
$\gamma_1< codZ-1$ and $Z$ is the only element in
$CLC(X,G_1)$.
Let $M\in |mH|$, for $m\gg 0$, a general member among
Cartier divisors containing $Z$. Let \hbox{$G^{\prime}=
G+\epsilon M$,}
and $\gamma^{\prime} =lct(X,G^{\prime})$. Define
$G_1:=1/\gamma^{\prime}
 G^{\prime}$ then $G_1\equiv \gamma_1 H$ and for $\epsilon\ll 1/m$ we
have $\gamma_1<cod Z$. Furthermore $(X,G_1)$ is LC and
$Z$ is the only element element of $CLC(X,G_1)$.
Since we have the strict inequality $\gamma_1<cod Z-1$, we may
furthermore assume,
by Kodaira Lemma, that there exists a log resolution $\mu:Y\ra X$
of $(X,G_1)$ such that
$$K_Y+E-A+F=\mu^*(K_X+G_1)+P,$$
where $\mu(E)=Z$, $A$ is an integral \hbox{$\mu$-exceptional}
divisor, $\lfloor F\rfloor=0$ and $P$ is an ample \Q-divisor.
Let
$$ N(t):= -K_Y-E-F+A+\mu^*(tH),$$
then $N(t)\equiv \mu^*(t+(n-2)-\gamma_1)H+P$ and $N(t)$ is ample
for
$t+(n-2)-\gamma_1> 0$, hence by hypothesis this is true whenever $t\geq
-n+1+cod Z$.
Thus K--V vanishing yields
\begin{equation}
H^i(Y,\mu^*(tH)-E+A)=0\hspace{.7cm}
H^i(E,(\mu^*(tH)+A)_{|E_0})=0
\label{van}
\end{equation}
for $i>0$ and $t\geq -n+1+cod Z$.
In particular there is the following surjection
$$
 H^0(Y,\mu^*H+A)\ra H^0(E,\mu^*H+A)\ra 0.
$$
Since $A$ is effective and $\mu$-exceptional, then
 $H^0(Y,\mu^*H+A)\simeq H^0(X,H)$, thus any section in $H^0(Y,\mu^*H+A)$
not vanishing on $E$, pushes forward to
give a section of $H$ not vanishing on $Z$.
To conclude the proof it is, therefore, enough to prove that
$h^0(E,\mu^*H+A)>0$. Let
$p(t)=\chi(E,\mu^*tH+A)$,
then by equation (\ref{van}),
$p(0)\geq 0$ and $p(t)=0$ for $0>t\geq -n+1+cod Z=-dim
Z+1$. Since $deg p(t)= dim Z$ and $p(t)>0$ for $t\gg 0$
then $h^0(E,\mu^*H+A)=p(1)>0$.
\eproof

We will first prove that log terminal Mukai varieties always have a log
terminal fundamental divisor.

\begin{Theorem} Let $X$ be a Mukai variety with log terminal
singularities and $K_X\equiv -(n-2)H$. Then the general element
in $|H|$ has log terminal singularities.
\label{muklt}
\end{Theorem}

\proof
By Proposition \ref{al} $dim |H|\geq 2$.
Let $S\in
|H|$ a generic
section and assume that $S$ has worse than LT singularities. Let
$\gamma=lct(X,S)$, then by our assumption
$\gamma\leq 1$, see for instance \cite[1.4]{Al}.
Let $Z\in CLC(X,\gamma S)$ a minimal center.
By Bertini Theorem $Z\subset Bsl|H|$ therefore by  Lemma \ref{KV}
$cod Z\leq 2$.

\claim{ $\gamma< cod Z$.}
\proof(of the Claim)
$\gamma\leq 1\leq cod Z$ and the equality could
hold only if $Z$ were a fixed
component of $|H|$
of multiplicity 1. Let	$S=Z+B$, then by connectedness
$W:=Z\cap B\neq
\emptyset$ and $W$ is properly contained in $Z$.
$S$ is a Cartier divisor singular along the codimension 2 subscheme
$W$ therefore $W\in CLC(X,S)$ and $Z$ cannot be minimal in $CLC(X,S)$.
\eproof
Let us perturb $S$, as in Lemma \ref{KV}, to construct a
\Q-divisor $S_1\equiv \gamma_1 H$ such that
$\gamma_1< codZ$ and $Z$ is the only element in
$CLC(X,\gamma S)$.

Let $\nu:Y\ra X$ a log resolution of $(X,S_1)$, with
$$K_Y-A+E+F=\nu^*(K_X+S_1),$$
where $A$ is $\nu$-exceptional, $\nu(E)=Z$ and $\lfloor F\rfloor=0$. Let
$$N:=\nu^*H+A-F-E-K_Y\equiv \nu^*(n-2)H,$$
then $N$ is nef and big and by K--V
vanishing we have,
\begin{equation}
H^1(Y,\nu^*H-E+A)=0 \label{van1}
\end{equation}

\claim{ $h^0(Z,H_{|Z})>0$}

Let us first prove that the claim gives us a contradiction.
Since $A$ does not contain $E$ and is effective then
$H^0(Z,H_{|Z})\hookrightarrow H^0(E,(\mu^*H+A)_{|E})$ therefore by the
vanishing (\ref{van1}) and the claim there exists a section
in $H^0(Y,\mu^*H+A)$
not vanishing on $E$.
Thus there exists a section of $H$ not vanishing on $Z$,
giving a contradiction and
proving the theorem.

\proof(of the claim)

$cod Z\leq 2$ thus we can apply subadjunction formula
of Theorem \ref{clc}.
There exists an  effective \Q-divisor
$\Delta$ on $Z$ such that
$$-(n-2-\gamma_1)H_{|Z}\equiv(K_X+S_1)_{|Z}\equiv K_Z+\Delta$$
and $(Z,\Delta)$ is KLT.
That is
to say
 $(Z,\Delta)$ is a log Fano variety with $i(Z)\geq dimZ-2$.
and the claim  follows directly from Proposition \ref{al}.
\eproof

The next step is to prove that canonical Mukai varieties have canonical
fundamental divisor.

\begin{Theorem} Let $X$ be a Mukai variety with canonical singularities.
Then the general element
in $|H|$ has canonical singularities.
\label{mukcan}
\end{Theorem}
\proof
By Theorem \ref{muklt} the general element $S\in |H|$ has LT singularities.
Let $\mu:Y\ra X$ a log resolution of $(X,S)$, with $\mu^*S=\overline{S}+
\sum r_i E_i$, where $\overline{S}$ is base point free,
 and $K_Y=\mu^*K_X+\sum a_i E_i$. Let us assume that $S$ has
not canonical singularities, then, maybe after reordering the indexes,
we have $a_0<r_0$. Since $S$ is generic then
$\mu(E_i)\subset Bsl|H|$, for all $i$ with $r_i>0$.
 Let $D=S+S^1$, with $S^1\in |H|$ a generic section.
First observe that $\mu$ is a log resolution of $(X,D)$.
Then $
(X,D)$ is not LC,
 infact $a_0+1<r_0+r^1_0$, where $r^1_0\geq 1$
 is the
multiplicity of $S_1$ at the center of the valuation associated to $E_0$.
Let $\gamma=lct(X,D)<1$ and $W$ a minimal center of
$CLC(X,\gamma D)$. $X$ has canonical singularities therefore whenever
$cod \mu(E_j)\leq 2$ then $a_j\in \N$.
 Thus, since $(X,D)$ is not LC and $S$ is LT, we have
$cod W\geq 3$.
We can therefore derive a contradiction by Lemma \ref{KV}.
\eproof

\Remark{ If $X$ is a terminal Mukai variety of dimension$\geq 4$
then by the above proposition
we immediately get that the generic element $S$
of the fundamental divisor $|H|$ is terminal.
In fact outside $Bsl|H|$
$S$ is terminal by Bertini Theorem and along the base locus the
discrepancy must be positive, since the generic section of
$H_{|S}$ is canonical.}

What remains to be done is to study terminal 3-folds with $-K_X\equiv H$,
let us start with some examples

\example Let us consider $X_{2,6}\subset \Proj(1,1,1,1,2,3)$ given
by the following equations
\begin{eqnarray*}
F_2(x_0,\ldots,x_3)=0\\
x_5^2+F_3(x_0,\ldots,x_4)x_5+F_6(x_0,\ldots,x_4)=0
\end{eqnarray*}
Assume that $F_6$ contains the monomial $x_4^3$ then
$$X\cap \{x_0=x_1=x_2=x_3=0\}=\{[(0:0:0:0:-1:1)]\}=\{p\},$$
thus $X$ is on the smooth locus
of $\Proj(1,1,1,1,2,3)$ and
$S$ and $Q$ are Cartier along $X$.
In particular adjunction formula holds $X$ is Gorenstein
and
$Bsl|-K_X|=\{p\}$.
Since $Q$ is singular at $p$ then $X$ is singular
at $p$ and therefore all elements in $|-K_X|$ are singular.

\example Let us now consider $Y_{2,4}\subset \Proj(1,1,1,1,1,2)$ given by
the following equations
\begin{eqnarray*}
x_5+F_2(x_0,\ldots,x_4)=0\\
F_4(x_0,\ldots,x_5)=0
\end{eqnarray*}
Since the first equation is linear then
$X$ is isomorphic to a quartic in $\Proj^4$. Let us choose this quartic
with two simple nodes at $(0,0,0,\pm1,1)$ and
consider the involution $\sigma$ on $\Proj(1,1,1,1,1,2)$ given by
$$(x_0:x_1:x_2:x_3:x_4:x_5)\mapsto(x_0:x_1:x_2:-x_3:-x_4:-x_5).$$
Let $\pi:Y\ra X=Y/ \sigma$ the quotient,
 then $X$ is a 3-fold with a $cA_1$ point,
\hbox{$p=\pi([(0:0:0:\pm1:1:0)])$,} and 8 points of singular index 2,
the fixed points of the involution. Furthermore $-K_X\equiv H$ and
$Bsl|H|=\{p\}$,
therefore all elements in $H$ are singular.

To prove the theorem we have to show that the
 above are the unique possibilities for a terminal Mukai variety which
has not good divisors.

\begin{Theorem} Let $X$ be a terminal Mukai 3-fold, assume that
all the divisor in the linear system $|H|$ are singular, then $X$ is one
of the following:
\begin{itemize}
\item[-]\cite{Mo}
 if $X$ is Gorenstein then $X$ is a complete intersection of a quadric
and a sestic in $\Proj(1,1,1,1,2,3)$
\item[-] if $X$ is not Gorenstein then the canonical cover of $X$ is a complete
intersection of a quadric and a quartic in $\Proj(1,1,1,1,1,2)$.
\end{itemize}
\end{Theorem}
\proof
Let $S\in |H|$ a generic element, then $S$ has at worse
canonical singularities, by Theorem \ref{mukcan}.
By K--V vanishing we get that $S$ is either a K3 surface or an
Enriques surface. Assume now that $Bsl|H|\neq\emptyset$ and
$S$ is singular, then by
\cite{SD} and \cite{Cos} we know that $S$ has an $A_1$ singularity
and $Bsl|H|=Bsl|H_{|S}|=\{x\}$, see also \cite{Shi} and \cite{Prok}. Let
$f:Y\ra X$ the blow up of the point $x$, with exceptional divisor $E$.
Then $f^*S=\tilde{S}+E$ and $\tilde{S}\cdot E^2=-2$ furthermore
$$0=f^*S\cdot E^2=\tilde{S}\cdot E^2+E^3,$$
thus $E^3=2$.  By \cite{SD} and
\cite{Cos}, $\tilde{S}_{|\tilde{S}}$ is an elliptic pencil, therefore
$\tilde{S}^3=0$. This gives
$$0=\tilde{S}^3=(f^*S-E)^3=\f^*H^3-2,$$
thus $H^3=2$, in case $X$ is not Gorenstein see also \cite{Prok}.

Assume that $X$ is not	Gorenstein then the generic element $E\in |H|$ is
a canonical Enriques surface with an
$A_1$ singularity. Let $\pi:Y\ra X$ the cyclic cover
associated to $\O_X(K_X+H)$, \cite{YPG}.
Let $S=\pi^{-1}E$ the pull back of a generic section $E\in |H|$
and $C=S_{|S}$. By connectedness $S$ is a canonical K3 surface
and $S^3=4$.
By Riemann--Roch theorem $h^0(C,S_{|C})=3$ and $h^0(C,nS_{|C})=2(2n-1)$,
for $n>1$.
Let $\{x_2,x_3,x_4\}$ be a basis of $H^0(C,S_{|C})$.

Let $\psi:S\flip S^{\prime}$ the map defined by the sections of
$H^0(S,C)$. If $\psi$ is not birational or is not a morphism, then,
by \cite[Th 5.2, Sect. 2.7]{SD} (see also \cite[Cor 2.2]{Shi}),
$H^0(C,S_{|C})^{\otimes 2}\neq
H^0(C,2S_{|C})$, in other words there is a  section
 $x_5\in H^0(C,2S_{|C})$ which is not in $H^0(C,S_{|C})^{\otimes 2}$, and
there is a quadratic relation of the kind $F_2(x_2,x_3,x_4,x_5)=0$.
 Going further we get that
nothing new happens
for $H^0(C,3S_{|C})$, while we get a relation in $H^0(C,4S_{|C})$,
of type \hbox{$F_4(x_2,x_3,x_4,x_5)=0$.}
This is enough to describe $Y$ as a complete
intersection of a quadric and a quartic in $\Proj(1,1,1,1,1,2)$. To
conclude observe that we can explicitly write down
an involution on $Y$ with only a finite number of fixed points as
$$(x_0:x_1:x_2:x_3:x_4:x_5)\mapsto (x_0:x_1:x_2:-x_3:-x_4:-x_5).$$

If the map $\psi$ is a birational
morphism then, by \cite[Th 6.1]{SD}, $Y$ is a
quartic in $\Proj^4$, embedded by $H^0(Y,\pi^*H)$. Note that
the involution $\pi$ must be
 the restriction of a projective transformation and
cannot, therefore, fix only finitely many points.

If $X$ is Gorenstein then the generic element $S\in |H|$
is a canonical K3 surface. Let $C=S_{|S}$
then with a similar calculation we get
that \hbox{$h^0(C,H_{|C})=2$} and $h^0(C,nH_{|C})=(2n-1)$, for $n>1$.
This time $H_{|S}$ is not base point free by  hypothesis, thus
there is a new section $x_4$
in $H^0(C,2H_{|C})$ which gives rise to a quadratic relation,
and a new
section $x_5$ in $H^0(C,3H_{|C})$, which gives rise to a relation in
$H^0(C,6H_{|C})$.
Therefore  we have the description
given in the proposition.
\eproof

What remains to be done is to prove Mukai hypothesis.
\begin{Theorem} Let $X$ be a smooth
Mukai
variety.
Then $X$ has good divisors.
\end{Theorem}

\proof If $X$ is smooth and the generic element in
$|H|$ is not smooth
then  by the surjection
$$H^0(X,H)\ra H^0(H,H_{|H})\ra 0,$$
and previous Theorems we know that $Bsl|H|=\{x\}$.
Let $H_i\in |H|$, for $i=1,\ldots, n-1$ generic elements and
$D=H_1+\cdots+H_{(n-1)}$, then the minimal center of $CLC(X,D)$ is $x$ and
$(X,D)$ is not LC at $x$, since $2(n-1)>n$.  We can therefore
derive a contradiction by Lemma \ref{KV}.
\eproof

 \begin{flushleft}
Massimiliano Mella\\
Universit{\`a} degli Studi di Trento\\
Dipartimento di Matematica\\
I38050 Povo (TN), Italia\\
e-mail: mella@science.unitn.it
\end{flushleft}
\end{document}